\documentclass{cmam}
\usepackage{amsmath,amssymb,amsfonts,xcolor}
\usepackage[ansinew]{inputenc}
\usepackage{epsfig,cite}
\usepackage {longtable}
\usepackage {comment}

\newcommand{{\un}}{u^{(0)}_{n}(x)}

\newcommand{{\bea}}{\begin{array}}
\newcommand{{\ea}}{\end{array}}
\newcommand{\be}{\begin{equation}}
\newcommand{\ee}{ \end{equation}}

\newcommand{{\il}}{\int}
\newcommand{{\la}}{\lambda_n^{(0)}}

\begin{document}
\runauthor{V.~Makarov and N.~Romaniuk} \runtitle{Symbolic Algorithm of the FD-method}
\begin{topmatter}

\title{Symbolic Algorithm of the Functional-Discrete Method for a Sturm-Liouville Problem with a Polynomial Potential}

\author{Volodymyr Makarov}
\address{Institute of Mathematics of NAS of Ukraine, 3 Tereshchenkivs’ka Str., 01004 Kyiv-4, Ukraine}
\email{makarov@imath.kiev.ua; http://orcid.org/0000-0002-4883-6574}
\author{Nataliia Romaniuk}
\address{Institute of Mathematics of NAS of Ukraine, 3 Tereshchenkivs’ka Str., 01004 Kyiv-4, Ukraine}
\email{romaniuknm@gmail.com; http://orcid.org/0000-0002-3497-7077}

\begin{abstract}
\noindent A new symbolic algorithmic implementation of the general scheme of the exponentially convergent functional-discrete (FD-) method is developed and justified for the Sturm-Liouville problem on a finite interval for the Schr{\"o}dinger equation with a polynomial potential and the boundary conditions of Dirichlet type. The algorithm of the general scheme of our method is developed when the potential function is approximated by the piecewise-constant function. Our algorithm is symbolic and operates with the decomposition coefficients of the eigenfunction corrections in some basis.
The number of summands in these decompositions depends on the degree of the potential polynomial and on the correction number.
Our method uses the algebraic operations only and does not need solutions of any boundary value problems and computations of any integrals unlike the previously version.  
The numerical example illustrates the theoretical results.

\end{abstract}
\date{August 11, 2017}
\end{topmatter}

\section{Introduction} \label{s0}

The eigenvalue problems (EVP) play an important role in various applications. There exist a large number of methods for their numerical solution.
Some of them are implemented 
in the application software packages included in such well-known algorithm and program 
libraries as NAG Fortran Library \cite{NAG_Mark20bib},  CALGO \cite{CALGObib}, CPC Program Library \cite{CPCbib}, JINRLIB 
\cite{JINRLIBbib}, etc. and in the collection of software packages of Ghent University \cite{GentUniversiteit}.

However, the many numerical methods possess disadvantages which are still only partially overcome.
They are the following: the accuracy degradation with the 
increasing of the eigenvalue index; usage of the mesh generated at the start of the 
numerical process; saturation of accuracy; the number of reliable numerical eigenvalues 
is limited and depends on a mesh step. 
Disadvantages of the classical discrete and spectral methods, including Numerov's method (NM), finite element (FEM) and finite difference methods (FDM), are outlined in the conclusion remarks in \cite{140_Zhang2015}
``Although the number of reliable eigenvalues increases with an increased computational scale \textit{N}, the percentage of reliable eigenvalues (compared with non-reliable eigenvalues) will go to zero when \textit{N} goes to infinity''. Here \textit{N} is the total degrees of freedom in resulting discrete 
systems.

The problem of the accuracy degradation with increasing of the eigenvalue index was 
partially overcome with the use of the asymptotic correction in conjunction with 
FDM, NM and FEM and was suggested by the Australian mathematicians J.~W.~Paine, F.~R.~de~Hoog, R.~S.~Anderssen and A.~L.~Andrew in the 80s of the 20th century. 
But, 
this approach is effective for the eigenvalues with very low indexes and isn't effective 
for the eigenvalues with large indexes (see the monograph \cite{GMR_Pryce1994} and the corresponding references therein). 
Meanwhile, there are  problems where a large number (thousands) of eigenvalues and eigenfunctions norms is required,  for example, for computing the spectral density function \cite[p.~273]{GMR_Pryce1994}.

For the last few years, the analytical (functional) methods based on the idea of the homotopy method (the parameter continuation method) (see e.g. \cite{118_GMR_AllgGeorg1990, 115_GMR_Arm1979}) are widely used for solving the eigenvalue problems. 
With the help for these methods one can find the solutions as fast convergent functional series. 
The major advantage of the analytical approximating methods is the possibility to study the properties of the solution of the original problem. Among these methods is 
the Adomian decomposition method (ADM) suggested in the 80s of the 20th century 
by the American physicist G.~Adomian (see references in \cite{65_Rach2012}). 
 Similar techniques are also applied in \cite{Liao_PhD_1992,30_He1999_VIM}.      
These methods have been applied to the problem of Sturm-Liouville type by S.~Abbasbandy (2011), M.~T.~Atay (2010), B.~S.~Attili and D.~Lesnic (2005--2007), S.~Irandoust-pakchin (2015), M.~A.~Jafari (2009), A.~Neamaty (2010), N.~Singh (2013), A.~H.~S.~Taher (2014) and others.

The development of the analytical methods in recent years has been associated with the improvement of the computer algebra systems, but nevertheless the computational costs of analytical methods can be quite high unlike the discrete methods. That is why the development of more efficient and less computationally costly algorithms for the high-precision analytical methods is an important and vital task today.

The idea of FD-method is closely related to above-mentioned methods 
HPM, HAM and ADM. FD-method in general sense has significant advantages over the 
purely analytical methods as it includes the discrete component with the help of 
which we can achieve the convergence when the mentioned methods HPM, HAM, ADM are 
divergent. 
Suggested by Volodymyr Makarov in 1991 the FD-method in \cite{MRL-Makarov1991} enables us to overcome the above-mentioned disadvantages 
of the discrete methods and can be applied to operator equations in general form. 
It was proved that in many cases the FD-method converges super-exponentially. The FD-method was developed by V.~L.~Makarov, I.~P.~Gavrilyuk, I.~I.~Lazurchak, N.~O.~Rossokhata, V.~B.~Vasylyk, D.~O.~Sytnyk, N.~M.~Romaniuk, D.~V.~Dragunov, A.~V.~Klymenko, B.~I.~Bandyrskii, O.~L.~Ukhanev, A.~M.~Popov and others.

The case when the potential function $q(x)$ is approximated by the constant zero $\bar{q}(x)\equiv 0$ (the simplest variant of the FD-method) can be considered 
as one of the variants of the homotopy method mentioned above \cite{118_GMR_AllgGeorg1990, 115_GMR_Arm1979}. 
The FD-method in general version 
(when the discrete component $\bar{q}(x)\not\equiv 0$) is also closely related to the 
methods which use the approximation to the coefficients of the differential equation. 
 The variants of these methods have been used since the beginning of the XX century, and for the piecewise-constant 
approximation for linear ordinary differential equations of second order it was justified 
(for the first time) in 1928 by N.~Bogoliouboff, N.~Kryloff in \cite{62_BogoliouboffKryloff1928} and was named the ``metodo dei tronconi''. 
Later 
in 1969 \cite{132_Gordon1969} the method of piecewise-polynomial 
approximation to the coefficients of the system of coupled linear second order differential 
equations was suggested by R.~G.~Gordon. 
The contribution of J.~D{\"a}hnn is also worth mentioning. In his article \cite{133_Dahnn1982}  a piecewise-constant approximation of the equation coefficients was used for the 
second order Sturm-Liouville problems. 
The methods of this type are known as ``Pruess 
methods'' \cite[Ch.6]{GMR_Pryce1994} and are named after S.~Pruess. 
He provided a thorough convergence and error analysis of such methods 
using a piecewise-polynomial approximation in 1973 \cite{117_GMR-Pruess1973}. 
The methods of the approximation to the coefficients of the differential 
equation also include exponentially weighted Legendre--Gauss Tau method (ELGT) \cite{35_DaouMatar2010,36_Daou2011} and piecewise perturbation methods \cite{GentUniversiteit}.

In the present paper the general scheme of the FD-method (when $\bar{q}(x)\not\equiv 0$) is applied to the Sturm-Liouville 
problem for the Schr{\"o}dinger differential 
equation
\begin{equation} \label{mr_1_} 
\frac{d^{2} u(x)}{dx^{2} } +(\lambda -q(x))
\, u(x)=0\, ,\, \, \; \; \; \; \; \; x\in (A,B)\, ,\, \; \; \; \; \; \; u(A)=u(B)=0, 
\end{equation} 
where $A,B$ are 
real constants, and the polynomial potential
\begin{equation} \label{mr_2_} 
q(x)=\sum _{p=0}^{r}c_{p} x^{p},\,\,\,\,c_{r}\ne 0. 
\end{equation}
In 
this case the FD-method for the problem \eqref{mr_1_}, \eqref{mr_2_} 
is exactly realizable in the sense that the corrections to the eigenfunctions $u_{n}^{(j+1)} 
(x)$ (see below Section~\ref{s1} below) can be expressed analytically in a closed form (see Definition in \cite{MRL-MakVinok1995}). 

A principally new symbolic algorithm of the simplest variant of the FD-method was proposed 
in \cite{MakarovRomanyuk2014} for a problem \eqref{mr_1_}, \eqref{mr_2_} and used in \cite{DemkGavrMak2016,GavrMakRom2017}
for several linear differential operators with fractional derivatives (for the case when the 
potential function $q(x)$ is approximated by the constant zero $\bar{q}(x)
\equiv 0$). The results below are the further development of these ideas.

Let us consider  the case when the simplest variant of the FD-method mentioned in 
\cite{MakarovRomanyuk2014} for finding the smallest 
eigenvalues of a problem \eqref{mr_1_}, \eqref{mr_2_} is divergent. In this case, the general scheme of the FD-method with a piecewise-constant function $\bar{q}(x)$ is used. 
 Our new algorithm
is a symbolic in this case in the sense that it uses only exact symbolic operations from a computer algebra system at each iteration step on the set of decomposition coefficients of the eigenfunction corrections in a certain basis.
Other the numerical values of target eigenvalue and the numerical values of the decomposition 
coefficients of the corresponding eigenfunction or their symbolic expressions can 
be considered the results of the algorithm execution. 
As a result FD-method could be implemented much more efficiently and provides the properties of "multivariance" and "multimodality" to the computational experiments (see \cite{Samarskij1979}).

Note, that the simplest variant of the FD-method for the problem \eqref{mr_1_}, 
\eqref{mr_2_} will be always convergent beginning with some eigenvalue index 
number (it can possibly be large enough). It means that  using the simplest variant 
of the FD-method one can obtain the asymptotic formulas for eigenvalues and eigenfunctions. 
Thus, for $q(x)=x$ we have the asymptotic formula
\[\lambda _{n} =(\pi n)^{2} +\frac{1}{2} +\frac{1}{48(\pi n)^{2} } -\frac{5}{16(
\pi n)^{4} } +\frac{1}{2304(\pi n)^{6} } -\frac{35}{384(\pi n)^{8} } +O(n^{-10} )\] 
In 
addition, the following formulas hold true for the corrections to any eigenvalue:
\begin{equation*} 
\begin{split}
&\lambda _{n}^{(2j)} =\frac{(4j-3)!}{2(3j)!(j-1)!3^{j-1} } \left(\frac{1}{2\pi n} 
\right)^{4j-2} -\frac{5(4j-1)!}{(3j-1)!(j-1)!3^{j} } \left(\frac{1}{2\pi n} \right)^{4j} 
+...,\\
&\lambda _{n}^{(1)} =\frac{1}{2} ,\; \; \; \lambda _{n}^{(2j+1)} =0,\; 
\; \; j=1,2,...,\; \; \; n=1,2,...
\end{split}
\end{equation*}

The paper is structured as follows. 
Section~\ref{s1} contains the traditional algorithm 
of the general scheme of the FD-method and the supplemental theoretical results needed 
for the developing of new symbolic algorithmic implementation mentioned in Section~\ref{s2}. Section IV describes the new approach. The theoretical results are illustrated by the given numerical example.

\section{Traditional algorithm of general scheme of FD-method} \label{s1}

We look for the exact solution of the eigenvalue problem \eqref{mr_1_}, \eqref{mr_2_} in the form of the series
\begin{equation} \label{mr_3_} 
u_{n} (x)=\sum _{j=0}^{\infty }u_{n}^{(j)}  (x),{
\kern 1pt} \; \; \; \lambda _{n} =\sum _{j=0}^{\infty }\lambda _{n}^{(j)}  ,\; \; 
\; n=1,2,... 
\end{equation} 
provided that they converge. The sufficient condition for the convergence 
of the series \eqref{mr_3_} will be presented below in Theorem~\ref{theor_mr_1}. We choose 
a mesh 
$$\omega =\left\{x_{i} ,i=1,2,...,N-1:\; A=x_{0} <\, x_{1} <...<\, x_{N} =B
\right\}.$$
The approximate solution to problem \eqref{mr_1_} is represented 
by a pair of corresponding truncated series
\begin{equation} \label{mr_4_} 
\mathop{u_{n}}\limits^{m}  (x)=\sum _{j=0}^{m}u_{n}^{(j)}  
(x),{\kern 1pt} \; \; \; \mathop{\lambda_{n}}\limits^{m}  =\sum _{j=0}^{m}\lambda 
_{n}^{(j)}  ,\; \; \; n=1,2,... 
\end{equation} 
which is called the approximation of rank $m$ \cite{MakarovRomanyuk2014,MRL-Makarov1991,Makarov1997,MRL-Band2000}. The summands of series \eqref{mr_3_}, 
\eqref{mr_4_} are the solutions of the recursive sequence of problems:
\begin{equation} \label{mr_5_} 
\begin{split}
&\frac{d^{2} u_{n}^{(j+1)} 
(x)}{dx^{2} } +(\lambda _{n}^{(0)} -\bar{q}( x))\, u_{n}^{(j+1)} (x)=F_{n}^{(j+1)} 
( x),\, \, x\in (A,B),\\
&u_{n}^{(j+1)} (A)=u_{n}^{(j+1)} (B)=0,\, \, \, j=-1,0,1,...,
\; \; \; n=1,2,...,
\end{split}
\end{equation} 
where
\begin{equation} \label{mr_6_} 
 F_{n}^{(j+1)} ( x)=-\sum _{p=0}^{j}\lambda _{n}^{(j+1-p)} \, u_{n}^{(p)} (x)+(q(
 \, x)-\bar{q}( x)) \, u_{n}^{(j)} (x),\; j=0,1,...,\,F_{n}^{(0)} 
 (x)=0.
\end{equation} 
Here the function $\bar{q}(x)$ is a piecewise-constant approximation to the potential 
polynomial function \eqref{mr_2_}, namely
\begin{equation} \label{mr_7_} 
\bar{q}(x)=\frac{q(x_{i-1} )+q( x_{i} 
)}{2} ,\, \, \, \; x\in \left[x_{i-1} ,x_{i} \right),\; \; i=\overline{1,N-1},\, 
\, \, x\in \left[x_{N-1} ,x_{N} \right].  
\end{equation} 
where the grid covers the closed interval $\left[A,B\right]$.
At the discontinuity points of $\bar{q}(x)$ the following matching conditions should be fulfilled:
\begin{equation} \label{mr_8_} 
\begin{split}
&\left. u_{n}^{(j+1)} (x)\right|_{x=x_{i} 
} =u_{n,i+1}^{(j+1)} (x_{i} )-u_{n,i}^{(j+1)} (x_{i} )=0,\\
&\left. \frac{du_{n}^{(j+1)} (x)}{dx} \right|_{x=x_{i} } =\left. \frac{du_{n,i+1}^{(j+1)} 
(x)}{dx} \right|_{x=x_{i} } -\left. \frac{du_{n,i}^{(j+1)} (x)}{dx} \right|_{x=x_{i} 
} =0,\\
&i=\overline{1,N-1},\, \, \, \; \; \; j=-1,0,1,...
\end{split}
\end{equation}

The start values $u_{n}^{(0)} (x)$, $\lambda _{n}^{(0)}$ for the recurrent process \eqref{mr_4_}--\eqref{mr_8_} one obtains from the basic problem \eqref{mr_5_}--\eqref{mr_8_} with $j=-1$. The function $u_{n}^{(0)} (x)$ possesses the representation
\begin{equation} \label{mr_9_} 
\begin{split}
&u_{n,i}^{(0)} (x)=\, \frac{a_{n,i,0}^{(0)} 
}{\kappa _{n,i} \, } \sin (\kappa _{n,i} \, x)+b_{n,i,0}^{(0)} \cos (\kappa _{n,i} 
\, x);\\
&x\in \left[x_{i-1} ,x_{i} \right)\, ,i=\overline{1,N-1};\; \; \; x\in 
\, \left[x_{N-1} ,x_{N} \right],i=N;\\
&\kappa _{n,i} =\sqrt{\lambda _{n}^{(0)} 
-\bar{q}(x_{i-1} )} \, ,\, \; \; \, i=\overline{1,N},\; \; \; n=1,2,...
\end{split}
\end{equation} 
with some constants  $a_{n,i,0}^{(0)}$, $b_{n,i,0}^{(0)}$, $i=\overline{1,N}$.
Substituting \eqref{mr_9_} into the boundary conditions in \eqref{mr_5_} and into the matching conditions \eqref{mr_8_} we obtain the following homogeneous system of the linear algebraic equations with a square matrix
$D(\lambda _{n}^{(0)})$ of order $2N$:
\begin{equation} \label{mr_10_} 
\left\{\begin{array}{l} {\, \frac{a_{n,1,0}^{(0)} 
}{\kappa _{n,1} \, } \sin (\kappa _{n,1} \, A)+b_{n,1,0}^{(0)} \cos (\kappa _{n,1} 
\, A)=0,} \\ {-\, \frac{a_{n,i,0}^{(0)} }{\kappa _{n,i} \, } \sin (\kappa _{n,i} 
\, x_{i} )-b_{n,i,0}^{(0)} \cos (\kappa _{n,i} \, x_{i} )+\frac{a_{n,i+1,0}^{(0)} 
}{\kappa _{n,i+1} \, } \sin (\kappa _{n,i+1} \, x_{i} )+} \\ {\;\;\;\;\;\;\;\;\;\;\;\;\;\;\;\;\;\;\;\;\;\;\;\;\;\;\;\;\;\;\;\;\;\;\;\;\;\;\;\;\;\;\;\;\;\;\;\;\;\;\;\;\;\;\;\;\;\;\;\;\;\;\;\;\;\;\;+b_{n,i+1,0}^{(0)} \cos (
\kappa _{n,i+1} \, x_{i} )=0,} \\ {-a_{n,i,0}^{(0)} \cos 
(\kappa _{n,i} \, x_{i} )+\kappa _{n,i} b_{n,i,0}^{(0)} \sin (\kappa _{n,i} \, x_{i} 
)+a_{n,i+1,0}^{(0)} \cos (\kappa _{n,i+1} \, x_{i} )-} \\ {\;\;\;\;\;\;\;\;\;\;\;\;\;\;\;\;\;\;\;\;\;\;\;\;\;\;\;\;\;\;\;\;\;\;\;\;\;\;\;\;\;\;\;\;\;\;\;\;\;\;\;\;\;\;\;\;\;\;\;\;\;\;\;\;\;\;\;-\kappa _{n,i+1} b_{n,i+1,0}^{(0)} 
\sin (\kappa _{n,i+1} \, x_{i} )=0,} \\ {\;\;\;\;\;\;\;\;\;\;\;\;\;\;\;\;\;\;\;\;\;\;\;\;\;\;\;\;\;\;\;\;\;\;\;\;\;\;\;\;\;\;\;\;\;\;\;\;\;\;\;\;\;\;\;\;\;\;\;\;\;\;\;\;\;\;\;\;\;\;\;\;\;\;\;\;\;\;\;\;\;\;\;\;\;\;\;\;\;\;\;\;\;\;\;\;i=\overline{1,N-1},} \\ {\frac{a_{n,N,0}^{(0)} 
}{\kappa _{n,N} \, } \sin (\kappa _{n,N}  B)+b_{n,N,0}^{(0)} \cos (\kappa _{n,N} 
\, B)=0.} \end{array}\right.  \end{equation} 
The homogeneous system of the linear equations \eqref{mr_10_} has a non-trivial 
solution if and only if its determinant is equal to zero, i.e.
\begin{equation} \label{mr_11_} 
\det (D(\lambda _{n}^{(0)} ))=0. 
\end{equation} 
We look for the roots of the equation \eqref{mr_11_} which are different from $\bar{q}(x_{i-1})$, $i=\overline{1,N}$ and each of which is the eigenvalue $\lambda_{n}^{(0)} $ of multiplicity 1 of the basic problem \eqref{mr_5_}--\eqref{mr_8_} with $j=-1$. 
For the given $\lambda _{n}^{(0)} $ the solution of the system \eqref{mr_10_} can be determined only up to a constant $b_{n,N,0}^{(0)} $ which we obtain from the normalization condition $\int _{A}^{B}[u_{n}^{(0)}(x)]^{2} dx=1 $ or, which is the same, $b_{n,N,0}^{(0)} $ is the root of the quadratic equation
\begin{equation} \label{mr_12_} 
\begin{split}
&\left(\cot ^{2} (\kappa _{n,N} 
)S_{n,N,0} +R_{n,N,0} \right)\left(b_{n,N,0}^{(0)} \right)^{2} +2 T_{n,N,0}\frac{ a_{n,N,0}^{(0)}}{\kappa _{n,N}} b_{n,N,0}^{(0)} =\\
&\; \; \; \; \; \; =1-\sum _{i=1}^{N-1} 
\left(\left(\frac{a_{n,i,0}^{(0)} }{\kappa _{n,i} \, } \right)^{2} S_{n,i,0} +2\frac{a_{n,i,0}^{(0)}}{\kappa _{n,i}}  b_{n,i,0}^{(0)} T_{n,i,0} +\left(b_{n,i,0}^{(0)} 
\right)^{2} R_{n,i,0} \right),
\end{split}
\end{equation} 
where we use the notations (see \cite{GradRyzh2015}): 
\begin{equation*} 
S_{n,i,p} =\int _{x_{i-1} }^{x_{i} }x^{p}  \sin ^{2} (\kappa 
_{n,i} x)dx=
\end{equation*}
\begin{equation*} 
=\frac{p!}{\kappa _{n,i}^{p+1} } \left(\frac{\left(\kappa _{n,i} x\right)^{p+1} 
}{2\left(p+1\right)!} +
\sum _{k=0}^{\left[\kern-0.15em\left[{p\mathord{\left/ {\vphantom 
{p 2}} \right. \kern-\nulldelimiterspace} 2} \right]\kern-0.15em\right]}\frac{(-1)^{k+1} 
\left(\kappa _{n,i} x\right)^{p-2k} }{(p-2k)!\, 2^{2k+2} }  \right. \sin (2\kappa 
_{n,i} x)+
\end{equation*}
\begin{equation}\label{mr_12_1} 
+\sum _{k=0}^{\left[\kern-0.15em\left[{(p-1)\mathord{\left/ {\vphantom 
{(p-1) 2}} \right. \kern-\nulldelimiterspace} 2} \right]\kern-0.15em\right]}\frac{(-1)^{k+1} 
\left(\kappa _{n,i} x\right)^{p-2k-1} }{(p-2k-1)!\, 2^{2k+3} }  \cos (2\kappa _{n,i} 
x)
-\left. \left. \frac{\left(-1\right)^{\left[\kern-0.15em\left[{(2p+1)\mathord{
\left/ {\vphantom {(2p+1) 4}} \right. \kern-\nulldelimiterspace} 4} \right]\kern-0.15em
\right]} \sin ^{2} \left({p\pi \mathord{\left/ {\vphantom {p\pi  2}} \right. \kern-
\nulldelimiterspace} 2} \right)}{2^{p+2} } \right)\right|_{x_{i-1} }^{x_{i} },
\end{equation}
\begin{equation*} 
\begin{split}
&R_{n,i,p} =\int _{x_{i-1} }^{x_{i} }x^{p}  \cos ^{2} (\kappa _{n,i} x)dx=\left. 
\frac{x^{p+1} }{p+1} \right|_{x_{i-1} }^{x_{i} } -S_{n,i,p} ,
\end{split}
\end{equation*}
\begin{equation*} 
\begin{split}
&T_{n,i,p} 
=\frac{1}{2} \int _{x_{i-1} }^{x_{i} }x^{p}  \sin (2\kappa _{n,i} x)dx=-\frac{p!}{2} 
\left. \sum _{k=0}^{p}\frac{x^{p-k} \cos \left(2\kappa _{n,i} x+{k\pi \mathord{\left/ 
{\vphantom {k\pi  2}} \right. \kern-\nulldelimiterspace} 2} \right)}{(p-k)!\, \left(2
\kappa _{n,i} \right)^{k+1} }  \right|_{x_{i-1} }^{x_{i} }.
\end{split}
\end{equation*}
In \eqref{mr_12_1} $p$ is an integer nonnegative number, $\left[\kern-0.15em\left[y\right]
\kern-0.15em\right]$ is the greatest integer less than or equal to a real number $y$.

Using the results of Chapter~12 in \cite{Atkinson1968} we conclude that the basic problem \eqref{mr_5_}--\eqref{mr_8_} with $j=-1$ possesses a countable set of the eigenvalues which are simple,  bounded below and with the infinity as the only accumulation point. 
Let us enumerate them in increasing order, i.e. 
$\lambda _{1}^{(0)} <\lambda _{2}^{(0)} <...<\lambda _{n}^{(0)}<...$. 
The sequence of the corresponding normalized eigenfunctions 
$\left\{u_{n}^{(0)} (x)\right\}_{n=1}^{\infty }$ builds a complete orthonormal system 
in Hilbert space $L_{2} \left[A,B\right]$. 

Using the completeness of the orthonormal system $\left\{u_{n}^{(0)} (x)\right\}_{n=1}^{
\infty}$ in $L_{2} \left[A,B\right]$ and the solvability condition 
\begin{equation} \label{mr_13_} 
\int _{A}^{B}F_{n}^{(j+1)}(x)u_{n}^{(0)}(x)dx=0,\; \; \; j=0,1,... 
\end{equation}
of problems \eqref{mr_5_}-\eqref{mr_8_} for fixed $j$  $(j=0,1,...)$, we obtain its solution in the form
\begin{equation} \label{mr_14_} 
u_{n}^{(j+1)} (x)=\sum _{p=1,\,p\ne n}^{\infty }\int _{A}^{B}F_{n}^{(j+1)} (\xi ) u_{p}^{(0)} (
\xi )d\xi  \frac{u_{p}^{(0)} (x)}{\lambda _{n}^{(0)} -\lambda _{p}^{(0)} } , 
\end{equation} 
\begin{equation}\label{mr_15_} 
\lambda _{n}^{(j+1)} =\int _{A}^{B}\left(q(x)-\bar{q}(x)\right)u_{n}^{(j)} (x)u_{n}^{(0)} (x)dx. 
\end{equation} 
The solution $u_{n}^{(j+1)} (x)$ satisfies the orthogonality condition 
\begin{equation}\label{mr_16_} 
\int _{A}^{B} u_{n}^{(0)} (x)u_{n}^{(j+1)} (x)dx=0,\; \; \; 
\; j=0,1,...
\end{equation}
Let us introduce the following notation
\begin{equation} \label{mr_17_} 
\left\| v\right\| _{\infty } =\mathop{\max}\limits_{x\in [A,B]} |v(x)|,\, \, \, \; \; \; \; \; \, \left\| v\right\| =\left
\{\int _{A}^{B}[v(x)]^{2} dx \right\}^{1/2} , 
\end{equation} 
\begin{equation} \label{mr_18} 
\begin{split}
&M_{n} =\max \left\{\left(\lambda _{n}^{(0)} 
-\lambda _{n-1}^{(0)} \right)^{-1} ,\left(\lambda _{n+1}^{(0)} -\lambda _{n}^{(0)} 
\right)^{-1} \right\},\, \, n=2,3,...,\\
&M_{1} =\left(\lambda _{2}^{(0)} -
\lambda _{1}^{(0)} \right)^{-1} .
\end{split}
\end{equation}
Then one can prove the following assertion.
\begin{theorem}{\cite{Makarov1997,MRL-Band2000}}\label{theor_mr_1}
Let the following condition holds true
\begin{equation}\label{mr_19_}
r_{n} 
=4\left\| q-\bar{q}\right\| _{\infty } M_{n} <1, n=1,2,...,                                          
\end{equation} 
then 
the FD-method for the Sturm-Liouville problem (1), (2) converges super-exponentially and  the 
following error estimates are valid:
\begin{equation} \label{mr_20_} 
\begin{split}
&\left|\lambda _{n} -\mathop{
\lambda _{n} }\limits^{m} \right|=\left|\lambda _{n} -\sum _{j=0}^{m}\lambda _{n}^{(j)}  
\right|\le \left\| q-\bar{q}\right\| _{\infty } \beta _{m} ,\\
&\left\| u_{n} -
\mathop{u_{n} }\limits^{m} \right\| \le \left\| u_{n} -\sum _{j=0}^{m}u_{n}^{(j)}  
\right\| \le \beta _{m+1},
\end{split}
\end{equation} 
\[m=0,1,2,...,\] 
where
\[\beta _{0} =1+\frac{r_{n} }{1-r_{n} } \alpha _{1} ,\; \; \; \; \; \beta _{m} =
\frac{r_{n} ^{m} }{1-r_{n} } \alpha _{m} ,\; \; m=1,2,...,\] 
\[\alpha _{m} =2\frac{(2m-1)!!}{(2m+2)!!} \le \frac{1}{(m+1)\sqrt{\pi m} } ,\; \; 
\; \; \; m=0,1,2,...\] 
\end{theorem}

\section{Symbolic algorithmic implementation of the FD-method} \label{s2}

According to \cite{MakarovRomanyuk2014}, for fixed $j$ $(j=-1,0,1,...,m-1)$ the  
solution of problem \eqref{mr_5_}--\eqref{mr_8_} is given by
\begin{equation} \label{mr_21_} 
u_{n,i}^{(j+1)} (x)=\, \sum _{p=0}^{M(j)}x^{p}  
\left(\frac{a_{n,i,p}^{(j+1)} }{\kappa _{n,i} \, } \sin (\kappa _{n,i} \, x)+b_{n,i,p}^{(j+1)} 
\cos (\kappa _{n,i} \, x)\right), 
\end{equation} 
\[x\in \left[x_{i-1} ,x_{i} \right)\, ,\; i=\overline{1,N-1};\; \; x\in \left[x_{N-1} 
,x_{N} \right],\, i=N;\] 
\[j=-1,0,1,...,m-1,\; \; \; n=1,2,...,\] 
where $M(j)=(j+1)(r+1)$, $a_{n,i,p}^{(j+1)}$ and $b_{n,i,p}^{(j+1)}$ $(p=\overline{0,M(j)})$ are 
the decomposition coefficients of the eigenfunction corrections $u_{n}^{(j+1)} (x)$ in the basis $x^{p} \sin (\kappa _{n,i} \, x)$, $x^{p} \cos (\kappa _{n,i} \, 
x)$ $(p=\overline{0,M(j)})$ on the each interval $\left[x_{i-1} ,x_{i} \right)\, 
,\; i=\overline{1,N-1},\; \left[x_{N-1} ,x_{N} \right]$ (see Definition 1). Unlike 
the \eqref{mr_14_}, the representation \eqref{mr_21_} is used below to  
develop a new symbolic algorithm of the FD-method. 

Using \eqref{mr_21_} in the boundary conditions of \eqref{mr_5_} 
and in the matching conditions \eqref{mr_8_} we obtain the following nonhomogeneous 
system of linear algebraic equations for the coefficients $a_{n,i,p}^{(j+1)} ,$ $b_{n,i,p}^{(j+1)}$, $p=\overline{0,M\left(j\right)}$, $i=\overline{1,N}$:
\begin{equation} \label{mr_23_} 
\left\{\!\begin{array}{l} {\frac{a_{n,1,0}^{(j+1)} 
}{\kappa _{n,1} } \sin (\kappa _{n,1}  A)+b_{n,1,0}^{(j+1)} \cos (\kappa _{n,1} 
A)=-\sum\limits _{p=1}^{M(j)}A^{p}  \left(\frac{a_{n,1,p}^{(j+1)} }{\kappa _{n,1} 
} \sin (\kappa _{n,1}  A)\right.+}\\
{\;\;\;\;\;\;\;\;\;\;\;\;\;\;\;\;\;\;\;\;\;\;\;\;\;\;\;\;\;\;\;\;\;\;\;\;\;\;\;\;\;\;\;\;\;\;\;\;\;\;\;\;\;\;\;\;\;\;\;\;\;\;\;\;\;\;\;\;\;\;\;\;\;\;\;\;\;\;\;\;\;\;\;\;\;\;\;\left.+b_{n,1,p}^{(j+1)} \cos (\kappa _{n,1}  A)\right),} 
\\ {-\frac{a_{n,i,0}^{(j+1)} }{\kappa _{n,i} } \sin (\kappa _{n,i}  x_{i} )-b_{n,i,0}^{(j+1)} 
\cos (\kappa _{n,i}  x_{i} )+\frac{a_{n,i+1,0}^{(j+1)} }{\kappa _{n,i+1} } \sin 
(\kappa _{n,i+1}  x_{i} )+}\\
{+b_{n,i+1,0}^{(j+1)} \cos (\kappa _{n,i+1}  x_{i} )=\sum _{p=1}^{M(j)}x_{i}^{p} \left(\frac{a_{n,i,p}^{(j+1)} }{\kappa _{n,i}  
} \sin (\kappa _{n,i}  x_{i} )+b_{n,i,p}^{(j+1)} \cos (\kappa _{n,i}  x_{i} )-
\right.}\\
{\;\;\;\;\;\;\;\;\;\;\;\;\;\;\;\;\;\;\;\;\;\;\;\;\;\;\;\;\;\;\;\;\;\;\;\;\;\;\;\;\;\;\;\;\;\;\;\;\;-\frac{a_{n,i+1,p}^{(j+1)} }{\kappa _{n,i+1}  } \sin (\kappa _{n,i+1} 
 x_{i} )\left.-b_{n,i+1,p}^{(j+1)} \cos (\kappa _{n,i+1} x_{i} )\right),} 
\\ {-a_{n,i,0}^{(j+1)} \cos (\kappa _{n,i} x_{i} )+b_{n,i,0}^{(j+1)} \kappa _{n,i} 
\sin (\kappa _{n,i} x_{i} )+a_{n,i+1,0}^{(j+1)} \cos (\kappa _{n,i+1}  x_{i} 
)-} \\ 
{\;\;\;\;\;\;\;\;\;\;\;\;\;\;\;\;-b_{n,i+1,0}^{(j+1)} \kappa _{n,i+1} \sin (\kappa _{n,i+1} x_{i} )=\sum 
_{p=1}^{M(j)}x_{i}^{p-1}  \left(a_{n,i,p}^{(j+1)} \left(\frac{p}{\kappa _{n,i}  
} \sin (\kappa _{n,i} x_{i} )\right.\right.+}\\
{\;\;\;\;\;\;\;\;\;\;\;\;\;\;\;\;\;\;\;\;\;\;\;\;\;\;\;\;\left.+x_{i} \cos (\kappa _{n,i} x_{i} )\right)+b_{n,i,p}^{(j+1)} 
\left(p\cos (\kappa _{n,i} x_{i} )-x_{i} \sin (\kappa _{n,i} x_{i} )\kappa _{n,i} 
\right)-}\\
{\;\;\;\;-a_{n,i+1,p}^{(j+1)} \left(\frac{p}{\kappa _{n,i+1} 
 } \sin (\kappa _{n,i+1} x_{i} )+x_{i} \cos (\kappa _{n,i+1} x_{i} )\right)-
\left. b_{n,i+1,p}^{(j+1)} 
\left(p\cos (\kappa _{n,i+1} x_{i} )-\right.\right.}\\
{\;\;\;\;\;\;\;\;\;\;\;\;\;\;\;\;\;\;\;\;\;\;\;\;\;\;\;\;\;\;\;\;\;\;\;\;\;\;\;\;\;\;\;\;\;\;\;\;\;\;\;\;\;\;\;\;\;\;\;\;\;\;\;\;\;\;\;\;\;\;\;\;\;\;\;\;\;\;\;\left.\left.-x_{i} \sin (\kappa _{n,i+1} x_{i} )\kappa _{n,i+1} 
\right)\right),}\\
{\;\;\;\;\;\;\;\;\;\;\;\;\;\;\;\;\;\;\;\;\;\;\;\;\;\;\;\;\;\;\;\;\;\;\;\;\;\;\;\;\;\;\;\;\;\;\;\;\;\;\;\;\;\;\;\;\;\;\;\;\;\;\;\;\;\;\;\;\;\;\;\;\;\;\;\;\;\;\;\;\;\;\;\;\;\;\;\;\;\;\;\;\;\;\;\;\;\;i=\overline{1,N-1},} \\ 
{\frac{a_{n,N,0}^{(j+1)} }{\kappa _{n,N}  
} \sin (\kappa _{n,N}  B)+b_{n,N,0}^{(j+1)} \cos (\kappa _{n,N} B)=-\sum _{p=1}^{M(j)}B^{p}  
\left(\frac{a_{n,N,p}^{(j+1)} }{\kappa _{n,N} } \sin (\kappa _{n,N} B)+\right.}\\
{\;\;\;\;\;\;\;\;\;\;\;\;\;\;\;\;\;\;\;\;\;\;\;\;\;\;\;\;\;\;\;\;\;\;\;\;\;\;\;\;\;\;\;\;\;\;\;\;\;\;\;\;\;\;\;\;\;\;\;\;\;\;\;\;\;\;\;\;\;\;\;\;\;\;\;\;\;\;\;\;\;\;\;\;\;\;\left.+b_{n,N,p}^{(j+1)} 
\cos (\kappa _{n,N} B)\right).} \end{array}\right.  
\end{equation}
The square matrix of this system is singular since it coincides with the matrix $D(
\lambda _{n}^{(0)} )$ of the system \eqref{mr_10_}. Let us introduce the following vectors
$$\vec{Y}_{n}^{(j+1)} =\left(a_{n,1,0}^{(j+1)} ,b_{n,1,0}^{(j+1)} 
,a_{n,2,0}^{(j+1)} ,b_{n,2,0}^{(j+1)} ,...,a_{n,N,0}^{(j+1)} ,b_{n,N,0}^{(j+1)} \right)^{T} 
,\;\;\;j=-1,0,1,...,m-1,$$
$$\vec{H}_{n}^{(j+1)} =\left(H_{n,1}^{(j+1)} ,H_{n,2}^{(j+1)} 
,H_{n,3}^{(j+1)} ,...,H_{n,2N-1}^{(j+1)} ,H_{n,2N}^{(j+1)} \right)^{T} ,\;\;\;j=0,1,...,m-1,
$$
$$\vec{H}_{n}^{(0)} =\vec{0}=(\underbrace{0,...,0}_{2N})^{T} ,$$
where
$$H_{n,1}^{(j+1)} =-\sum _{p=1}^{M(j)}A^{p}  \left(\frac{a_{n,1,p}^{(j+1)} }{\kappa 
_{n,1} \, } \sin (\kappa _{n,1} \, A)+b_{n,1,p}^{(j+1)} \cos (\kappa _{n,1} \, A)
\right),$$
\begin{equation*}
\begin{split}
&H_{n,2i}^{(j+1)} =\sum _{p=1}^{M(j)}x_{i}^{p} \left(\frac{a_{n,i,p}^{(j+1)} 
}{\kappa _{n,i} \, } \sin (\kappa _{n,i} \, x_{i} )+b_{n,i,p}^{(j+1)} \cos (\kappa 
_{n,i} \, x_{i} )-\right.\\
&\; \; \; \; \; \; \; \; \; \; \; \; -\left. \frac{a_{n,i+1,p}^{(j+1)} 
}{\kappa _{n,i+1} \, } \sin (\kappa _{n,i+1} \, x_{i} )-b_{n,i+1,p}^{(j+1)} \cos 
(\kappa _{n,i+1} \, x_{i} )\right),
\end{split}
\end{equation*} 
\begin{equation*}
\begin{split}
&H_{n,2i+1}^{(j+1)} =\sum_{p=1}^{M(j)}x_{i}^{p-1}  \left(a_{n,i,p}^{(j+1)} 
\left(\frac{p}{\kappa _{n,i} \, } \sin (\kappa _{n,i} x_{i} )+x_{i} \cos (\kappa 
_{n,i} x_{i} )\kappa _{n,i} \right)+\right.\\
&\; \; \; \; \; \; \; \; \; \; \; 
\; \; \; \; \; \; \; \; \; \; \; \; +b_{n,i,p}^{(j+1)} \left(p\cos (\kappa _{n,i} 
x_{i} )-x_{i} \sin (\kappa _{n,i} x_{i} )\kappa _{n,i} \right)-\\
&\; \; \; \; 
\; \; \; \; \; \; \; \; \; \; \; \; \; \; \; \; \; -a_{n,i+1,p}^{(j+1)} \left(\frac{p}{
\kappa _{n,i+1} \, } \sin (\kappa _{n,i+1} x_{i} )+x_{i} \cos (\kappa _{n,i+1} x_{i} 
)\kappa _{n,i+1} \right)-\\
&\; \; \; \; \; \; \; \; \; \; \; \; \; \; \; \; \; 
\; \; \; \; \; \; \; \left. -b_{n,i+1,p}^{(j+1)} \left(p\cos (\kappa _{n,i+1} x_{i} 
)-x_{i} \sin (\kappa _{n,i+1} x_{i} )\kappa _{n,i+1} \right)\right),
\end{split}
\end{equation*}
$$i=\overline{1,N-1},$$
$$H_{n,2N}^{(j+1)} =-\sum _{p=1}^{M(j)}B^{p}  \left(\frac{a_{n,N,p}^{(j+1)} }{
\kappa _{n,N} \, } \sin (\kappa _{n,N} \, B)+b_{n,N,p}^{(j+1)} \cos (\kappa _{n,N} 
\, B)\right).$$
Now the systems \eqref{mr_10_} and \eqref{mr_23_} can 
be presented in the following matrix-vector form
\begin{equation} \label{mr_24_} 
D(\lambda _{n}^{(0)} )\vec{Y}_{n}^{(j+1)} =
\vec{H}_{n}^{(j+1)} ,\; \; j=-1,0,1,...,m-1,\; \; \; \vec{H}_{n}^{(0)} =\vec{0}. \end{equation} 
Because 
the eigenvalues of the basic problem \eqref{mr_5_}--\eqref{mr_8_} with $j=-1$ as well as the eigenvalue of the matrix $D(\lambda _{n}^{(0)} )$ with the index zero is simple, 
the necessary and sufficient condition for the solvability of \eqref{mr_24_} 
is the orthogonality of its right-hand side $\vec{H}_{n}^{(j+1)}$ to the eigenvector $\vec{Z}_{n} $ of 
the matrix of the corresponding conjugate homogeneous system \cite{Gantmacher1998}:
\begin{equation} \label{mr_25_} 
\vec{Z}_{n}^{T} \vec{H}_{n}^{(j+1)} =0,\; \; 
j=0,1,...,m-1. 
\end{equation} 
It can be shown that the condition \eqref{mr_25_} is equivalent to condition 
\eqref{mr_13_}.

Using the formulas \eqref{mr_21_} and \eqref{mr_15_} (it follows from 
(13) or \eqref{mr_25_}) we obtain the following formula for the corrections 
of the eigenvalues:
\begin{equation} \label{mr_27_} 
\begin{split}
&\lambda _{n}^{(j+1)} =\sum 
_{i=1}^{N} \sum _{p=0}^{r}c_{p} \sum _{t=0}^{M(j-1)}  \left(\frac{a_{n,i,t}^{(j)} 
a_{n,i,0}^{(0)} }{(\kappa _{n,i} )^{2} \, } \left(S_{n,i,t+p} -S_{n,i,t} \frac{x_{i-1}^{p} 
+x_{i}^{p} }{2} \right)\right. +\\
&\; \; \; \; \; \; \; \; \; \; \; \; \; \; +
\frac{a_{n,i,t}^{(j)} b_{n,i,0}^{(0)} +b_{n,i,t}^{(j)} a_{n,i,0}^{(0)} }{\kappa _{n,i} 
} \left(T_{n,i,t+p} -T_{n,i,t} \frac{x_{i-1}^{p} +x_{i}^{p} }{2} \right)+\\
&\; 
\; \; \; \; \; \; \; \; \; \; \; \; \; +b_{n,i,t}^{(j)} b_{n,i,0}^{(0)} \left. \left(R_{n,i,t+p} 
-R_{n,i,t} \frac{x_{i-1}^{p} +x_{i}^{p} }{2} \right)\right).
\end{split}
\end{equation} 

Let 
us find the formulas for the decomposition coefficients of representation \eqref{mr_21_}. 
For this purpose we substitute \eqref{mr_21_} into \eqref{mr_5_}, 
\eqref{mr_6_} and change the summation order in \eqref{mr_6_}. The result 
is
\begin{equation} \label{mr_28_} 
F_{n,i}^{(j+1)} (x)=\sum _{p=0}^{M(j)-1} x^{p} 
\left(f_{n,i,p}^{(j+1)} \sin (\kappa _{n,i} x)+g_{n,i,p}^{(j+1)} \cos (\kappa _{n,i} 
x)\right), i=\overline{1,N}, 
\end{equation} 
\begin{equation} \label{mr_29_} f_{n,i,p}^{(j+1)} =\sum _{l=\max (0,p-M(j-1))}^{
\min (r,p)}c_{l}  \frac{a_{n,i,p-l}^{(j)} }{\kappa _{n,i} \, } ,\; \; \; \; \; \; 
\; \; \; \; g_{n,i,p}^{(j+1)} =\sum _{l=\max (0,p-M(j-1))}^{\min (r,p)}c_{l} b_{n,i,p-l}^{(j)}  
, \end{equation} 
\[p=\overline{M(j)-r,M(j)-1}\, ,\] 
\begin{equation} \label{mr_30_} 
\begin{split}
&f_{n,i,t}^{(j+1)} =\sum _{l=
\max (0,t-M(j-1))}^{\min (r,t)}c_{l}  \, \frac{a_{n,i,t-l}^{(j)} }{\kappa _{n,i} 
\, } -\sum _{s=\left. \left. \right]\kern-0.12em\right]\frac{t}{r+1} \left[\kern-0.12em
\left[\right. \right. }^{j}\lambda _{n}^{(j+1-s)}  \frac{a_{n,i,t}^{(s)} }{\kappa 
_{n,i} \, } -\frac{q(\, x_{i-1} )+q(\, x_{i} )}{2} \frac{a_{n,i,t}^{(j)} }{\kappa 
_{n,i} \, },\\
&g_{n,i,t}^{(j+1)} =\sum _{l=\max (0,t-M(j-1))}^{\min (r,t)}c_{l}  
b_{n,i,t-l}^{(j)} \, -\sum _{s=\left. \left. \right]\kern-0.12em\right]\frac{t}{r+1} 
\left[\kern-0.12em\left[\right. \right. }^{j}\lambda _{n}^{(j+1-s)}  b_{n,i,t}^{(s)} 
-\frac{q(\, x_{i-1} )+q(\, x_{i} )}{2} b_{n,i,t}^{(j)} ,
\end{split}
\end{equation}
\[t=
\overline{0,M(j)-r-1},\] 
where $\left. \left. \right]\kern-0.15em\right]y\left[\kern-0.15em\left[\right. \right. $ is 
the smallest integer greater than or equal to a real number $y$. We require that the polynomials in the front of corresponding trigonometric functions are equal 
on the both sides of equation \eqref{mr_5_}, \eqref{mr_28_}. 
It leads to the recurrence system for the unknown coefficients of representation 
\eqref{mr_21_}:
\begin{equation} \label{mr_31_} 
\begin{split}
&(t+1)\left((t+2)\frac{a_{n,i,t+2}^{(j+1)} 
}{\kappa _{n,i} \, } -2\kappa _{n,i} b_{n,i,t+1}^{(j+1)} \right)=f_{n,i,t}^{(j+1)} 
,\\
&(t+1)\left((t+2)b_{n,i,t+2}^{(j+1)} +2a_{n,i,t+1}^{(j+1)} \right)=g_{n,i,t}^{(j+1)} 
,\\
&t=\overline{0,M(j)-2},\\
&-2M(j)\kappa _{n,i} b_{n,i,M(j)}^{(j+1)} =f_{n,i,M(j)-1}^{(j+1)} 
=c_{r} \frac{a_{n,i,M(j-1)}^{(j)} }{\kappa _{n,i} \, } ,\\
&2M(j)a_{n,i,M(j)}^{(j+1)} 
=g_{n,i,M(j)-1}^{(j+1)} =c_{r} b_{n,i,M(j-1)}^{(j)} .
\end{split}
\end{equation}
The consequence of this system is the following

\begin{equation} \label{mr_32_} 
\begin{split}
&(t+1)(t+2)\frac{a_{n,i,t+2}^{(j+1)} 
}{\kappa _{n,i} \, } +4\kappa _{n,i} a_{n,i,t}^{(j+1)} =f_{n,i,t}^{(j+1)} +\frac{2
\kappa _{n,i} }{t} g_{n,i,t-1}^{(j+1)} ,\\
&(t+1)(t+2)b_{n,i,t+2}^{(j+1)} +4\kappa 
_{n,i} ^{2} b_{n,i,t}^{(j+1)} =g_{n,i,t}^{(j+1)} -\frac{2\kappa _{n,i} }{t} f_{n,i,t-1}^{(j+1)} 
,\\
&t=\overline{1,M(j)-2},\\
&-2M(j)\kappa _{n,i} b_{n,i,M(j)}^{(j+1)} =f_{n,i,M(j)-1}^{(j+1)} 
=c_{r} \frac{a_{n,i,M(j-1)}^{(j)} }{\kappa _{n,i} \, } ,\\
&2M(j) \frac{a_{n,i,M(j)}^{(j+1)}}{\kappa _{n,i} }
=g_{n,i,M(j)-1}^{(j+1)} =c_{r} b_{n,i,M(j-1)}^{(j)} ,
\end{split}
\end{equation} 
Introducing the new variables by
\begin{equation} \label{mr_33_} 
\frac{a_{n,i,M(j)-2p}^{(j+1)} }{\kappa _{n,i} 
\, } =v_{p}^{a},\;\;\;\;  b_{n,i,M(j)-2p}^{(j+1)} =v_{p}^{b}, 
\end{equation} 
we obtain the following two linear Cauchy problems \textbf{(with $t=a$ and $t=b$)} for the linear inhomogeneous first order difference equation
\begin{equation} \label{mr_34_} 
v_{p-1}^{t} +\frac{4\kappa _{n,i} ^{2} }{(M(j)-2p+1)(M(j)-2p+2)} 
v_{p}^{t} =L_{p}^{t} ,\; \; \; t=a,b, p=\overline{1,\left[\kern-0.15em\left[{M(j)
\mathord{\left/ {\vphantom {M(j) 2}} \right. \kern-\nulldelimiterspace} 2} \right]
\kern-0.15em\right]-1} 
\end{equation} 
with the initial conditions
\begin{equation} \label{mr_35_} 
v_{0}^{a} =\frac{a_{n,i,M(j)}^{(j+1)} }{\kappa 
_{n,i} \, } =\frac{g_{n,i,M(j)-1}^{(j+1)} }{2M(j)\kappa _{n,i} } =\frac{c_{r} b_{n,i,M(j-1)}^{(j)} 
}{2M(j)\kappa _{n,i} } ,
\end{equation} 
\begin{equation} \label{mr_35_1} 
v_{0}^{b} =b_{n,i,M(j)}^{(j+1)} =-\frac{f_{n,i,M(j)-1}^{(j+1)} 
}{2M(j)\kappa _{n,i} } =-\frac{c_{r} a_{n,i,M(j-1)}^{(j)} }{2M(j)(\kappa _{n,i} )^{2} 
}  
\end{equation} 
and with the nonhomogeneities
\begin{equation} \label{mr_36_} L_{p}^{a} =\frac{1}{(M(j)-2p+1)(M(j)-2p+2)} 
\left(f_{n,i,M(j)-2p}^{(j+1)} +\frac{2\kappa _{n,i} }{M(j)-2p} g_{n,i,M(j)-2p-1}^{(j+1)} 
\right), \end{equation} 
\begin{equation} \label{mr_37_} L_{p}^{b} =\frac{1}{(M(j)-2p+1)(M(j)-2p+2)} 
\left(g_{n,i,M(j)-2p}^{(j+1)} -\frac{2\kappa _{n,i} }{M(j)-2p} f_{n,i,M(j)-2p-1}^{(j+1)} 
\right). \end{equation} 
The solution of \eqref{mr_34_}, \eqref{mr_35_}, \eqref{mr_36_} (with $t=a$) and \eqref{mr_34_}, 
\eqref{mr_35_1}, \eqref{mr_37_} (with $t=b$) is the expression
\begin{equation}\label{mr_38_} 
v_{p}^{t} =\frac{(-1)^{p} }{(2\kappa _{n,i} )^{2p} (M(j)-2p)!} \left(v_{0}^{t} 
M(j)!-\sum _{s=0}^{p-1} (-1)^{s} (2\kappa _{n,i} )^{2s} (M(j)-2s)!L_{s+1}^{t} \right),
\end{equation} 
\[t=a,b,\;\;\;\;\;p=\overline{1,\left[\kern-0.15em\left[{M(j)\mathord{
\left/ {\vphantom {M(j) 2}} \right. \kern-\nulldelimiterspace} 2} \right]\kern-0.15em
\right]-1}. \] 

Substituting  in \eqref{mr_31_}
\begin{equation} \label{mr_39_} 
\frac{a_{n,i,M(j)-2p+1}^{(j+1)} }{\kappa _{n,i} 
\, } =w_{p}^{a},\;\;\;\;\;   b_{n,i,M(j)-2p+1}^{(j+1)} =w_{p}^{b},
 \end{equation} 
we obtain the following two initial value problems (with $t=a$ and $t=b$) for the  linear inhomogeneous first-order difference equation
\begin{equation} \label{mr_40_} 
w_{p-1}^{t} +\frac{4(\kappa _{n,i} )^{2} }{(M(j)-2p+2)(M(j)-2p+3)} 
w_{p}^{t} =L_{p-{1\mathord{\left/ {\vphantom {1 2}} \right. \kern-\nulldelimiterspace} 
2} }^{t} ,\; \; \; t=a,b,\;\;\;p=\overline{1,\left[\kern-0.15em\left[{M(j)\mathord{
\left/ {\vphantom {M(j) 2}} \right. \kern-\nulldelimiterspace} 2} \right]\kern-0.15em
\right]} 
\end{equation} 
with the nonhomogeneities $L_{p-{1\mathord{\left/ {\vphantom {1 2}} 
\right. \kern-\nulldelimiterspace} 2} }^{a} ,\; \; L_{p-{1\mathord{\left/ {\vphantom 
{1 2}} \right. \kern-\nulldelimiterspace} 2} }^{b} $ (evaluated by substituting $p-{1
\mathord{\left/ {\vphantom {1 2}} \right. \kern-\nulldelimiterspace} 2}$ instead of $p$ into \eqref{mr_36_}, 
\eqref{mr_37_}) and with the corresponding initial conditions
\begin{equation} \label{mr_41_} 
w_{0}^{a} 
=\frac{a_{n,i,M(j)+1}^{(j+1)} }{\kappa _{n,i} \, } =0,
\end{equation} 
\begin{equation} \label{mr_41_1} 
w_{0}^{b} =b_{n,i,M(j)+1}^{(j+1)} 
=0. 
\end{equation} 
The solution of \eqref{mr_40_}, \eqref{mr_41_} (with $t=a$) and \eqref{mr_40_}, \eqref{mr_41_1} (with $t=b$) is the expression:
\begin{equation} \label{mr_42_}  
w_{p}^{t} =\frac{(-1)^{p+1} }{(2\kappa _{n,i} )^{2p} (M(j)-2p+1)!} \sum _{s=0}^{p-1} 
(-1)^{s} (2\kappa _{n,i} )^{2s} (M(j)-2s+1)!L_{s+{1\mathord{\left/ {\vphantom {1 
2}} \right. \kern-\nulldelimiterspace} 2} }^{t} ,
\end{equation} 
$$t=a,b,\;\;\;\; p=\overline{1,\left[\kern-0.15em\left[{M(j)
\mathord{\left/ {\vphantom {M(j) 2}} \right. \kern-\nulldelimiterspace} 2} \right]
\kern-0.15em\right]}.$$
Returning to the replacements \eqref{mr_33_}, \eqref{mr_39_} and the 
notations \eqref{mr_36_}, \eqref{mr_37_}, we obtain  the following recursive 
representation for the coefficients in \eqref{mr_21_} (see the notations \eqref{mr_29_}, 
\eqref{mr_30_}):
\begin{equation} \label{mr_43_} 
a_{n,i,M(j)}^{(j+1)} =\frac{c_{r} b_{n,i,M(j-1)}^{(j)} 
}{2M(j)} ,\;\;\;\;b_{n,i,M(j)}^{(j+1)} =-\frac{c_{r} a_{n,i,M(j-1)}^{(j)} }{2M(j)(\kappa 
_{n,i} )^{2} } , 
\end{equation}
\begin{equation} \label{mr_44_} 
\begin{split}
a_{n,i,M(j)-2p}^{(j+1)}& =
\frac{(-1)^{p} }{(2\kappa _{n,i} )^{2p} (M(j)-2p)!} \left(\frac{c_{r} b_{n,i,M(j-1)}^{(j)} 
(M(j)-1)!}{2} \right. -\\
&-\left. \kappa _{n,i} \sum _{s=0}^{p-1}(-1)^{s}  
(2\kappa _{n,i} )^{2s} (M(j)-2s-3)!\left((M(j)-2s-2)f_{n,i,M(j)-2s-2}^{(j+1)} +\right.\right.\\
&\left.+2\kappa _{n,i} g_{n,i,M(j)-2s-3}^{(j+1)} \right)\Bigg),
\end{split}
\end{equation} 
\begin{equation} \label{mr_45_} 
\begin{split}
b_{n,i,M(j)-2p}^{(j+1)}& =\frac{(-1)^{p+1} 
}{(2\kappa _{n,i} )^{2p} (M(j)-2p)!} \left(\frac{c_{r} a_{n,i,M(j-1)}^{(j)} (M(j)-1)!}{2(
\kappa _{n,i} )^{2} } \right. +\\
&+\left. \sum 
_{s=0}^{p-1}(-1)^{s}  (2\kappa _{n,i} )^{2s} (M(j)-2s-3)!\left((M(j)-2s-2)g_{n,i,M(j)-2s-2}^{(j+1)}-\right.\right.\\
&\left.-2\kappa _{n,i} f_{n,i,M(j)-2s-3}^{(j+1)} \right)\Bigg),
\end{split}
\end{equation} 
\[p=
\overline{1,\left[\kern-0.15em\left[{M(j)\mathord{\left/ {\vphantom {M(j) 2}} \right. 
\kern-\nulldelimiterspace} 2} \right]\kern-0.15em\right]-1},\]
\begin{equation} \label{mr_46_} 
\begin{split}
a_{n,i,M(j)-2t+1}^{(j+1)}& 
=\frac{(-1)^{t+1} \kappa _{n,i} }{(2\kappa _{n,i} )^{2t} (M(j)-2t+1)!} \times\\
&\times \sum _{s=0}^{t-1}(-1)^{s}  (2\kappa 
_{n,i} )^{2s} (M(j)-2s-2)!\left((M(j)-2s-1)f_{n,i,M(j)-2s-1}^{(j+1)}+\right.\\
&\left.+2\kappa _{n,i} 
g_{n,i,M(j)-2s-2}^{(j+1)} \right),
\end{split}
\end{equation} 
\begin{equation} \label{mr_47_} 
\begin{split}
b_{n,i,M(j)-2t+1}^{(j+1)}& 
=\frac{(-1)^{t+1} }{(2\kappa _{n,i} )^{2t} (M(j)-2t+1)!} \times\\
&\times \sum _{s=0}^{t-1}(-1)^{s}  (2\kappa _{n,i} )^{2s} 
(M(j)-2s-2)!\left((M(j)-2s-1)g_{n,i,M(j)-2s-1}^{(j+1)}-\right.\\
&\left. -2\kappa _{n,i} f_{n,i,M(j)-2s-2}^{(j+1)} 
\right),
\end{split}
\end{equation} 
\[t=\overline{1,\left[\kern-0.15em\left[{M(j)\mathord{\left/ {\vphantom {M(j) 2}} 
\right. \kern-\nulldelimiterspace} 2} \right]\kern-0.15em\right]}.\]

The coefficients $a_{n,i,0}^{(j+1)} ,b_{n,i,0}^{(j+1)} ,i=\overline{1,N}$ are determined 
by the system of the linear nonhomogeneous algebraic equations \eqref{mr_23_}, 
where the right-hand sides can be evaluated according to \eqref{mr_43_}--
\eqref{mr_47_}. As it was mentioned above, the matrix of the system \eqref{mr_23_} 
is singular, since it coincides with the matrix $D(\lambda _{n}^{(0)} )$ of the 
system \eqref{mr_10_}. 

The solution of the system \eqref{mr_21_} can be determined only up to the 
constant $b_{n,N,0}^{(0)} $, which we calculate from the orthogonality 
condition \eqref{mr_16_}, which leads to the following equation
\begin{equation} \label{mr_48_} 
\sum _{i=1}^{N} \sum _{t=0}^{M(j)}\left(\left(\frac{a_{n,i,0}^{(0)} }{\kappa _{n,i} 
} S_{n,i,t} +b_{n,i,0}^{(0)} T_{n,i,t} \right)\frac{a_{n,i,t}^{(j+1)} 
}{\kappa _{n,i} } +\left(\frac{a_{n,i,0}^{(0)} }{\kappa _{n,i} } T_{n,i,t} +b_{n,i,0}^{
(0)} R_{n,i,t} \right)b_{n,i,t}^{(j+1)} \right) =0.
\end{equation} 
with respect to the unknown $b_{n,N,0}^{(0)}$.

Thus, 
the formulas \eqref{mr_9_}-\eqref{mr_12_}, \eqref{mr_21_}, \eqref{mr_23_}, 
\eqref{mr_27_} and \eqref{mr_43_}-\eqref{mr_48_} represent the 
new symbolic algorithmic implementation of the general scheme of the FD-method for 
the problem \eqref{mr_1_} with the polynomial potential \eqref{mr_2_}. 
Our method uses the algebraic operations only and does not need the solutions 
of any boundary value problems \eqref{mr_5_}-\eqref{mr_8_} and computations of any integrals \eqref{mr_14_}-\eqref{mr_17_} unlike the previously known traditional implementations of the FD-method  \cite{MRL-Makarov1991,Makarov1997,MRL-Band2000}. 
This new scheme also improves a symbolic algorithm of the simplest variant of the FD-method in the case  $\bar{q}(x)\equiv 0$, $N=1$, proposed earlier in \cite{MakarovRomanyuk2014}.


\section{Numerical example}\label{s3}

\begin{example} 
We consider the eigenvalue problem \eqref{mr_1_} with $A=0,$ $B=1$ and with 
the potential \eqref{mr_2_}, where $r=1$, $c_{0} =-60,$ $c_{1} =120$ ($q(x)=-60+120x$). 
The computations of the exact eigenvalues (further denoted by $\lambda _{n}^{ex} $) 
and the approximations (denoted by $\mathop{\lambda_{n}}\limits^{m}  $) have been 
done with the help of the computer algebra system Maple 17 (Digits=200). The first 
four smallest exact eigenvalues of the problem under consideration are the following
$$\lambda 
_{1}^{ex} = -3.08815211843854844862886684381,$$
$$\lambda _{2}^{ex} = 41.5266775137315677830945919694,$$
$$\lambda _{3}^{ex} = 91.4591579961161898490753991651,$$
$$\lambda _{4}^{ex} = 159.625216916146830891863813793.$$

For the numerical results obtained by the FD-method of rank $m=\overline{0,10},\; 
20$ see Tables~1 and 2. Here the absolute errors 
\begin{equation}\label{mr_50_} 
\Delta _{n}^{(m)} =\; 
|\mathop{\lambda _{n} }\limits^{m} -\lambda _{n}^{ex} |
\end{equation}
of the approximations  $\mathop{
\lambda _{n} }\limits^{m} $ to the exact eigenvalues $\lambda _{n}^{ex} $ and the 
norms of the corresponding residuals are given by
\begin{equation}\label{mr_51_} 
\Omega _{n}^{(m)} =\left\| \Phi _{n}^{(m)} (x)\right\| =\left\{\int _{0}^{1}[\Phi 
_{n}^{(m)} (x)]^{2} dx \right\}^{1/2}
\end{equation}
with 
\[\Phi _{n}^{(m)} (x)=\frac{d^{2} \mathop{u_{n}}\limits^{m}  (x)}{dx^{2} } +\left(
\mathop{\lambda _{n} }\limits^{m} +60-120x\right)\mathop{u_{n}}\limits^{m}  (x).\] 

According 
to the Theorem~\ref{theor_mr_1} if the potential function $q(x)$ is approximated by $\bar{q}(x)\equiv 0$ ($N=1$) the sufficient convergence condition \eqref{mr_19_} is fulfilled for the  eigenpairs with the index $n\ge 13$. For $n=\overline{1,12}$ the FD-method can be divergent. However, as can be seen 
in Table~1, the simplest variant of the FD-method converges beginning with $n
\ge 3$ but for $n=1,2$ the FD-method is divergent. It means that the conditions of 
Theorem~\ref{theor_mr_1} are rough and can be improved.

Table~2 contains the results of the new symbolic algorithm described in Section~\ref{s2}. Here the potential function $q(x)$ was approximated by the piecewise constant approximation $\bar{q}(x)$ in the following two ways: 
I) the interval $(0,1)$ was partitioned into two equal subintervals ($N=2$, $x_{1} =\frac{1}{2}$); 
II) the interval $(0,1)$ was partitioned into three equal subintervals ($N=3$, $x_{1} =\frac{1}{3}$, $x_{2} =\frac{2}{3}$). 
The numerical results for the first two eigenpairs with $n=1,2$ are given in Table~2. One can observe that the convergence rate increases with increases together with increasing of the number of subdivision points (from one to two) and increases together with the index of the eigenpair  (see Table~2).

\begin{table}
\begin{center}
\begin{tabular}{|p{0.1in}|p{0.6in}|p{0.6in}|p{0.6in}|p{0.6in}|p{0.6in}|p{0.6in}|p{0.6in}|p{0.6in}|} \hline 
$m$ & \multicolumn{2}{|p{1.2in}|}{$n=1$} & \multicolumn{2}{|p{1.2in}|}{$n=2$} & \multicolumn{2}{|p{1.2in}|}{$n=3$} & \multicolumn{2}{|p{1.2in}|}{$n=4$} \\ \hline 
 & $\Delta 
_{1}^{(m)} $ & $\Omega _{1}^{(m)} $ & $\Delta _{2}^{(m)} $ & $\Omega _{2}^{(m)} $ & $\Delta 
_{3}^{(m)} $ & $\Omega _{3}^{(m)} $  & $\Delta _{4}^{(m)} $ & $\Omega _{4}^{(m)} $  \\ \hline 
0 & 13.0 & 21.7 & 2.05 & 31.9 & 2.63 & 33.5 & 1.71 & 34.0 \\ \hline 
1 & 13.0 & 23.5 & 2.05 & 11.4 & 2.63 & 12.5 & 1.71 & 9.86 \\ \hline 
2 & 2.84 & 9.00 & 2.66 & 5.87 & 0.174 & 3.03 & 7.77e-3 & 2.33 \\ \hline 
3 & 2.84 & 9.59 & 2.66 & 5.41 & 0.174 & 1.05 & 7.77e-3 & 5.96e-1 \\ \hline 
4 & 1.76 & 4.87 & 1.80 & 3.72 & 3.14e-2 & 2.74e-1 & 2.05e-3 & 1.36e-1 \\ \hline 
5 & 1.76 & 6.72 & 1.80 & 3.54 & 3.14e-2 & 1.10e-1 & 2.05e-3 & 3.7e-2 \\ \hline 
6 & 1.50 & 4.21 & 1.51 & 3.25 & 1.11e-2 & 3.67e-2 & 2.85e-5 & 8.43e-3 \\ \hline 
7 & 1.50 & 6.35 & 1.51 & 3.27 & 1.11e-2 & 2.22e-2 & 2.85e-5 & 2.34e-3 \\ \hline 
8 & 1.48 & 3.87 & 1.48 & 3.06 & 1.81e-3 & 4.52e-3 & 1.37e-5 & 5.45e-4 \\ \hline 
9 & 1.48 & 6.40 & 1.48 & 3.33 & 1.81e-3 & 3.75e-3 & 1.37e-5 & 1.50e-4 \\ \hline 
10 & 1.59 & 4.35 & 1.59 & 3.48 & 7.10e-5 & 4.26e-4 & 7.51e-7 & 3.46e-5 \\ \hline 
20 & 4.16 & 10.0 & 4.16 & 8.23 & 5.29e-7 & 4.42e-7 & 1.21e-12 & 2.93e-11 \\ \hline 
\end{tabular}
\caption{Numerical results obtained by the simplest variant of the FD-method ($\bar{q}(x)\equiv 0,$ $N=1$) of rank $m=\overline{0,10},20$: absolute errors $\Delta _{n}^{(m)}$ and norms of the corresponding residuals $\Omega _{n}^{(m)}$ (see notation \eqref{mr_50_}, \eqref{mr_51_}) for the eigenpairs $\lambda _{n},u_{n} (x),n=\overline{1,4}$.}
\end{center}
\end{table}

\begin{table}
\begin{center}
\begin{tabular}{|p{0.1in}|p{0.6in}|p{0.6in}|p{0.6in}|p{0.6in}|p{0.6in}|p{0.6in}|p{0.6in}|p{0.6in}|} \hline 
$m$ & \multicolumn{4}{|p{1.2in}|}{$n=1$} & \multicolumn{4}{|p{1.2in}|}{$n=2$} \\ \hline 
 & \multicolumn{2}{|p{1.2in}|}{I) $N=2,$ \newline $x_{1} =\frac{1}{2}$ } & \multicolumn{2}{|p{1.2in}|}{II) $N=3$, \newline $x_{1} =\frac{1}{3}$, $x_{2} =\frac{2}{3}$} & \multicolumn{2}{|p{1.2in}|}{I) $N=2,$ \newline $x_{1} =\frac{1}{2}$} & \multicolumn{2}{|p{1.2in}|}{II) $N=3$, \newline $x_{1} =\frac{1}{3}$, $x_{2} =\frac{2}{3}$} \\ \hline 
 & $\Delta 
_{1}^{(m)} $ & $\Omega _{1}^{(m)} $  & $\Delta _{1}^{(m)} $ & $\Omega _{1}^{(m)} $ & $\Delta 
_{2}^{(m)} $ & $\Omega _{2}^{(m)} $  & $\Delta _{2}^{(m)} $ & $\Omega _{2}^{(m)} $ \\ \hline 
0 & 2.77 & 15.2 & 1.11 & 11.3 & 6.61 & 15.8 & 1.20 & 10.7 \\ \hline 
1 & 1.20 & 1.74 & 0.278 & 0.375 & 7.27e-1 & 1.65 & 3.01e-1 & 4.61e-1 \\ \hline 
2 & 3.64e-2 & 0.115 & 2.52e-3 & 1.33e-2 & 5.39e-2 & 0.312 & 2.13e-3 & 1.43e-2 \\ \hline 
3 & 6.28e-3 & 1.31e-2 & 1.29e-4 & 5.35e-4 & 1.02e-2 & 5.51e-2 & 3.63e-4 & 1.16e-3 \\ \hline 
4 & 8.26e-4 & 1.81e-3 & 4.44e-6 & 2.03e-5 & 1.68e-3 & 6.98e-3 & 9.71e-6 & 6.81e-5 \\ \hline 
5 & 4.41e-5 & 2.13e-4 & 2.78e-8 & 1.17e-6 & 1.35e-4 & 6.53e-4 & 1.22e-6 & 3.99e-6 \\ \hline 
6 & 1.81e-5 & 3.27e-5 & 2.46e-9 & 3.38e-8 & 2.32e-5 & 3.44e-4 & 7.22e-8 & 2.57e-7 \\ \hline 
7 & 4.28e-7 & 4.51e-6 & 3.54e-11 & 2.21e-9 & 1.47e-5 & 1.13e-4 & 3.03e-9 & 1.33e-8 \\ \hline 
8 & 3.55e-7 & 7.10e-7 & 1.21e-11 & 1.38e-10 & 4.31e-6 & 2.49e-5 & 3.28e-10 & 8.86e-10 \\ \hline 
9 & 3.62e-8 & 1.04e-7 & 7.19e-13 & 5.55e-12 & 8.19e-7 & 3.45e-6 & 1.04e-11 & 5.06e-11 \\ \hline 
10 & 5.50e-9 & 1.69e-8 & 2.35e-14 & 4.44e-13 & 7.00e-8 & 5.04e-7 & 1.79e-12 & 4.29e-12 \\ \hline 
20 & 1.08e-16 & 2.20e-16 & 2.23e-26 & 1.96e-25 & 1.74e-14 & 7.81e-13 & 6.36e-24 & 1.67e-23 \\ \hline 
\end{tabular}
\caption{Numerical results obtained by the new symbolic algorithm of the general scheme of the FD-method of rank $m=\overline{0,10},20$ with the piecewise constant approximation $\bar{q}(x)$ in cases: 
      I) $N=2,$ $x_{1} =\frac{1}{2}$; 
      II) $N=3$, $x_{1} =\frac{1}{3}$, $x_{2} =\frac{2}{3}$ 
   for the eigenpairs $\lambda _{n},$ $u_{n} (x),$ $n=1,2$ (see notations \eqref{mr_50_}, \eqref{mr_51_}).}
\end{center}
\end{table}

\end{example} 

\newpage

\bibliographystyle{elsarticle-num}

\bibliography{SAFD2608}{}

\end{document}